\documentclass[11pt]{amsart}
\usepackage{amsmath}
\usepackage{amssymb}
\usepackage{graphicx}
\usepackage{cases}
\usepackage{fancyhdr}
\usepackage[inner=3cm,outer=4.5cm]{geometry}
\newtheorem{theorem}{Theorem}[section]

\theoremstyle{definition}

\numberwithin{equation}{section}
\makeatother
\title[A characterization of an integral triangle inequality]{\vspace{2cm}A characterization of a nonlinear integral triangle inequality}
\author[Ahmed A. Abdelhakim]{Ahmed A. Abdelhakim}
\address[]{
\hspace{-0.4cm}{Ahmed A. Abdelhakim}\newline
Department of Mathematics, College of Sciences and Arts in Unaizah, Qassim University, Qassim, Saudi Arabia\newline Mathematics Department, Faculty of Science, Assiut University, Assiut 71516, Egypt}
\email{{a.mouhamed@qu.edu.sa,ahmed.abdelhakim@aun.edu.eg}}
\begin{document}
\setcounter{page}{1165}
\begin{abstract}
Let $(E,\|.\|)$ be a Banach space and let $(\Omega,\mu)$ be a Lebesgue measure space. We characterize, for all $p>0$, measurable functions $u:\Omega\rightarrow \mathbb{R}$ for which
\begin{equation*}
\left\|
\int_{\Omega} f\,d\mu
\right\|^{p}\,\leq\,\int_{\Omega} u \|
 f \|^{p}\,d\mu.\tag{I}
\end{equation*}
We characterize $u$ for the reverse of (I) as well. The discrete counterpart of this problem is also solved.
\end{abstract}
\maketitle
\thispagestyle{empty}
\section{Introduction}
Let $p>0$. Suppose that $(E,\|.\|)$ is a Banach space and $(\Omega,\mu)$ is a Lebesgue measure space.
We characterize real-valued  measurable functions
$u$ on $\Omega$ such that the generalized integral triangle inequality
\begin{equation}\label{m1}
\left\|
\int_{\Omega} f\,d\mu
\right\|^{p}\,\leq\,\int_{\Omega} u \|
 f \|^{p}\,d\mu
\end{equation}
holds for all $f\in L^{p}(\Omega\rightarrow E)$.
We also characterize $u:\Omega\rightarrow \mathbb{R}$  for which the reverse of (\ref{m1}) holds.
For $p> 1$, the sufficient conditions on $u$ are obtained by manipulating (\ref{m1}) and applying the reverse H\"{o}lder's inequality.
The necessity of these conditions, on the other hand, result from considering an isoperimetric functional minimization problem. The remaining cases are treated using mostly counterexamples.  \\
\indent A particular discrete analogue
that motivates this problem was handled recently
by Dadipour et al. \cite{Dadipour}.
Their arguments are based on a finite dimensional generalization of
the results in \cite{Takagi} that provide a characterization for the envelope of prescribed families of lines in $\mathbb{R}^{2}$.
The results in \cite{Dadipour} can be summarized as follows:
\begin{theorem}\label{Dadi}
Let $(X,\|.\|)$ be a normed space, $x_{i}\in X,$ $i=1,...,n$ and $\lambda:=(\lambda_{1},...,\lambda_{n})\in \mathbb{R}^{n}$. Then $\,\|\sum_{i=1}^{n} x_{i}\|^{p}
\leq \sum_{i=1}^{n}{{\lambda_{i}} \| x_{i}\|^{p}}\,$
if and only if
\begin{numcases}{}
\lambda_{i}>0,\;\sum_{i=1}^{n}
\lambda_{i}^{\frac{1}{1-p}}\leq 1, & $p>1$;\label{pgt11}\\
\lambda_{i}\geq 1, & $0<p\leq1$.\label{pleq11}
\end{numcases}
While the reverse $\,\|\sum_{i=1}^{n} x_{i}\|^{p}
\geq \sum_{i=1}^{n}{{\lambda_{i}}\| x_{i}\|^{p}}\,$
holds if and only if $\lambda_{i}<0$ or
\begin{numcases}{\exists! j:\lambda_{j}>0,\:
\lambda_{i}<0,\,i\neq j,\;}
\lambda^{\frac{1}{1-p}}_{j}\leq 1 +\sum_{
i\neq j}|\lambda_{i}|^{\frac{1}{1-p}}, & $p>1$;\label{pgt1}\\
\lambda_{j}\leq
\min\{1,\,|\lambda_{i}|,i\neq j\}, & $0<p\leq1$.\label{pleq1}
\end{numcases}
\end{theorem}
We can certainly borrow, from (\ref{pgt11}) and (\ref{pleq11}) of Theorem \ref{Dadi},
the necessary conditions on the real sequence $\alpha:=\left(\alpha_{i}\right)_{i\geq 1}$ for the generalized triangle inequality
\begin{equation}\label{m21}
\|\sum_{i\geq 1} x_{i}\|^{p}
\leq \sum_{i\geq 1}{\alpha_{i} \| x_{i}\|^{p}},\;
x_{i} \in E,\;p> 0,
  \end{equation}
to hold. Nevertheless, the approach in \cite{Dadipour} is not helpful
in deriving sufficient conditions on $\alpha$ for (\ref{m21}) to hold.
When $p>1$, we benefit from the validity of  H\"{o}lder's inequality in any measure space, and identify the desired sufficient conditions on $\alpha$, in a treatment similar to that of (\ref{m1}). These conditions coincide, it turns out, with the necessary conditions suggested by Theorem \ref{Dadi}.
We also characterize $\alpha$ for
which (\ref{m21}) holds when $0<p\leq 1$, and
$\alpha$ for
which the reverse of (\ref{m21}) holds when $p>0$
as well.\\
\indent
The difference in nature between the Lebesgue
and counting measures shows in the difference between the characteristics of $u$ of (\ref{m1})
and those of $\alpha$ in (\ref{m21}).
\section{The integral inequality}
Denote by $\mu$ the Lebesgue measure on
$\mathbb{R}^{n}$ and let $\Omega$ be $\mu$-measurable. Let $(E,\|.\|)$ be a Banach space and consider the Banach space
$L^{1}(\Omega\rightarrow E)$ of strongly (Bochner) measurable functions $f:\Omega\rightarrow E$ that are also Bochner integrable (see \cite{Yosida}). Every $f\in L^{1}(\Omega\rightarrow E)$ satisfies
\begin{equation}\label{boch1}
\left\|
\int_{\Omega} f\,d\mu
\right\|\,\leq\,\int_{\Omega} \|
 f \|\,d\mu\,<\infty.
\end{equation}
In each of the remaining arguments in this section, $v$ will denote an arbitrary but prescribed nonzero vector in $E$. Suppose that $u:\Omega\rightarrow \mathbb{R}$ is $\mu$-measurable.
First, let's get the easier case $p=1$ out of the way:
\begin{theorem}
\begin{equation*}
\left\|
\int_{\Omega} f\,d\mu
\right\|\,\leq\,\int_{\Omega} u\|
 f \|\,d\mu,\; \forall\, f\in L^{1}(\Omega\rightarrow E),
\end{equation*}
if and only if
$u\geq 1$ $\mu$-a.e.$\,\Omega$.
\end{theorem}
\begin{proof}
The assertion follows readily from
the sharpness of the integral triangle inequality (\ref{boch1}). To establish necessity, suppose there exists $\,\widetilde{\Omega}\subset \Omega\,$ such that $\,\mu(\widetilde{\Omega})>0$ and $u<1$ $\mu$-a.e. $\widetilde{\Omega}$, then the piecewise constant $E$-valued function
\begin{equation*}
f(x):=\left\{
        \begin{array}{ll}
          v, & \hbox{$x\in\widetilde{\Omega}$;} \\
          0, & \hbox{$x\in\Omega\setminus \widetilde{\Omega}$}
        \end{array}
      \right.
\end{equation*}
would satisfy
$\,\displaystyle
\int_{\Omega} f\,d\mu\,=\,\frac{v}{\|v\|}
\int_{\widetilde{\Omega}} \|v\|\,d\mu
\,=\ \frac{v}{\|v\|}
\int_{\widetilde{\Omega}} \|f\|\,d\mu,$
and the contradiction
\begin{equation*}
\left\| \int_{\Omega} f\,d\mu \right\|
\,> \,
\int_{\widetilde{\Omega}}u\|f\|\,d\mu
\,= \,\int_{{\Omega}}u\|f\|\,d\mu
\end{equation*}
would follow.
\end{proof}
Secondly, we characterize $u$ for the reverse of
(\ref{m1}) when $p=1$. Actually, for any $p>0$, there does not exist $\mu$-almost everywhere positive functions $u$ for which the opposite of the inequality (\ref{m1}) always holds.
\begin{theorem}Let $p>0$. Then
\begin{equation}\label{oppo}
\left\|
\int_{\Omega} f\,d\mu
\right\|^{p}\,\geq\,\int_{\Omega} u\|
 f \|^{p}\,d\mu,\; \forall\, f\in L^{1}(\Omega\rightarrow E),
\end{equation}
if and only if
\begin{equation}\label{necrev}
u\leq 0,\; \mu\text{-a.e.}\,\Omega.
\end{equation}
\end{theorem}
\begin{proof}
Suppose that (\ref{oppo}) holds while
$\Omega$ has a subset ${\Omega}_{0}$
with $\mu({\Omega}_{0})>0$ and such that $u>0$
on ${\Omega}_{0}$.
Let $f$ be such that
\begin{equation}\label{ex}
\int_{\Omega} f\,d\mu\neq 0,\;
\frac{\|
 f \|}{\|
 \int_{\Omega} f\,d\mu \|}\,\geq \,\frac{c({\Omega}_{0})}{\epsilon}\,\chi_{{\Omega}_{0}},\; \epsilon>0,
\end{equation}
where $c({\Omega}_{0})$ is a constant that may depend on ${\Omega}_{0}$  but is independent of $\epsilon$.
We have to choose $\epsilon$ small enough for (\ref{boch1}) to be satisfied. Then, the inequality (\ref{oppo}) necessitates that
$\,\int_{{\Omega}_{0}}u\,d\mu\leq \epsilon^{p}/c^{p}({\Omega}_{0}).$ Therefore, by the arbitrary smallness
of $\epsilon$, we have $\;\int_{{\Omega}_{0}}u\,d\mu\,=\,0,\,$
a contradiction. \\
\indent For an example of a function that satisfies
(\ref{ex}), let $\,{\Omega}_{0}={\Omega}_{1}\cup{\Omega}_{2},\,$
where $\,{\Omega}_{1}\cap{\Omega}_{2}=\emptyset,\,$
and $\,\mu({\Omega}_{1})=\mu(
{\Omega}_{2}).$ Define
\begin{equation*}
 g(x):=\left\{
   \begin{array}{ll}
    (1+\epsilon)v , & \hbox{$x\in{\Omega}_{1}$;} \\
-v, & \hbox{$x\in{\Omega}_{2}$;} \\
     0, & \hbox{$x\in\Omega\setminus{\Omega}_{0}$.}
   \end{array}
 \right.
\end{equation*}
One can see that $\,\|g\|=\|v\|
\left((1+\epsilon)\chi_{{\Omega}_{1}}+\chi_{{\Omega}_{2}} \right),\,$ and
\begin{equation*}
  \int_{\Omega}g\,d\mu \,=\,
(1+\epsilon) \mu({\Omega}_{1})v
-\mu({\Omega}_{2})v\,=\,\epsilon \mu({\Omega}_{1})v.
\end{equation*}
Thus $g$ satisfies (\ref{ex}) with $c({\Omega}_{0})=2/\mu({\Omega}_{0}).$\\
\indent A different counterexample that shows the necessity of (\ref{necrev}) comes from considering
\begin{equation*}
h(x):=\left\{
        \begin{array}{ll}
          v, & \hbox{$x\in{\Omega}_{1}$;} \\
          -v, & \hbox{$x\in{\Omega}_{2}$;} \\
          0, & \hbox{$x\in\Omega\setminus{\Omega}_{0}$,}
        \end{array}
      \right.
\end{equation*}
where the $\Omega_{i}$'s are as before.
We have
$\,\|h\|=\|v\| \chi_{{\Omega}_{0}}\,$ and
$\,\int_{\Omega}h\,d\mu\,=\,0$, which, assuming (\ref{oppo}), leads to the contradiction $\|\int_{\Omega}h\,d\mu\|^{p}=0\geq\|v\|^{p} \int_{{\Omega}_{0}} u\,d\mu.$
\end{proof}
For $p\neq 1$, we have the subsequent characterization:
\begin{theorem}\label{main1}
Let $p>1$. Then
\begin{equation}\label{pg11}
\left\|
\int_{\Omega} f\,d\mu
\right\|^{p}\,\leq\,\int_{\Omega} u \|
 f \|^{p}\,d\mu,\;\;\forall\, f\in L^{1}(\Omega\rightarrow E)
\end{equation}
if and only if
\begin{equation}\label{assum1}
u>0\;\mu\text{-a.e.}\,\Omega,\;\; \int_{\Omega}\,u^{\frac{1}{1-p}}d\mu\leq 1.
\end{equation}
\end{theorem}
\begin{proof}
We may assume that $f\neq 0$ $\mu$-almost everywhere in $\Omega$. Let $u>0\,$ $\mu$\text{-a.e.}$\,\Omega$.
If $\int_{\Omega} f\,d\mu= 0$, then (\ref{pg11})
is fulfilled. Now let $\int_{\Omega} f\,d\mu\neq 0$. Invoking the reverse H\"{o}lder's inequality yields
\begin{eqnarray*}
\int_{\Omega} u \left(\frac{\| f \|}{\left\|
\int_{\Omega} f d\mu
\right\|}\right)^{p}
d\mu\geq
\left(\int_{\Omega} \frac{\| f \|}{\left\|
\int_{\Omega} f d\mu
\right\|}
\,d\mu\right)^{p}\left(\int_{\Omega} u^{\frac{1}{1-p}}\,d\mu\right)^{1-p}\geq1
\end{eqnarray*}
when $\,\int_{\Omega}\,u^{\frac{1}{1-p}}d\mu\leq 1,\,$ by (\ref{boch1}).
Hence (\ref{assum1}) is sufficient for (\ref{pg11}).\\
\indent Conversely, assume that (\ref{pg11}) holds true. Arguing by contradiction, let $\bar{\Omega}\subset \Omega$, $\mu(\bar{\Omega})>0$ and
$u\leq 0$ on $\bar{\Omega}$. Let $f:= v \chi_{\bar{\Omega}}$. Then
\begin{equation*}
\left\|\int_{\Omega} f\,d\mu \right\|^{p} =\|v\|^{p} \left(\mu(\bar{\Omega})\right)^{p},\;\;
\int_{\Omega}u \|f\|^{p}\,d\mu=
\|v\|^{p}\int_{\bar{\Omega}}u \,d\mu.
\end{equation*}
And (\ref{pg11}) implies
$0\geq\int_{\bar{\Omega}}u \,d\mu\geq \left(\mu(\bar{\Omega})\right)^{p}>0$.
Therefore, strict positivity of $u$
$\mu$-almost everywhere in $\Omega$
is necessary for (\ref{pg11}).\\
\indent Let $f(x):=\varphi(x)\,v$
where $\varphi\in L^{1}(\Omega\rightarrow \mathbb{R})$
has constant sign $\mu$-a.e.$\,\Omega$, and define
$\psi(x):=\varphi(x)/\int_{\Omega}\varphi\,d\mu$.
Observe that $\psi>0$ $\mu$-a.e.$\,\Omega$
and that $\int_{\Omega}\psi\,d\mu=1$.
If (\ref{pg11}) holds then
\begin{equation}\label{fn}
\int_{\Omega} u \left(\frac{\| f \|}{\left\|
\int_{\Omega} f d\mu
\right\|}\right)^{p}
d\mu=\int_{\Omega} u \frac{|\varphi|^{p}}
{|\int_{\Omega}\varphi d\mu|^{p}}\,d\mu=\int_{\Omega} u\, \psi^{p}\,d\mu
\geq 1.
\end{equation}
Let $\mathcal{A}:=\{\psi:\Omega\rightarrow \mathbb{R}:
\psi>0,\:\mu\text{-a.e.}\,\Omega,\;
\int_{\Omega}\psi\,d\mu=1\}$. Consider the functional
\begin{equation}\label{functional}
J[\psi]:=\int_{\Omega} u\psi^{p}\,d\mu,\;\psi \in
\mathcal{A}.
\end{equation}
In the light of (\ref{fn}), the inequality (\ref{pg11})
implies
\begin{equation}\label{JJ}
 J[\psi]\geq 1,\;\;\forall\psi\in \mathcal{A}.
\end{equation}
Let
\begin{equation*}
\psi_{0}:=\frac{u^{\frac{1}{1-p}}}{\int_{\Omega}
u^{\frac{1}{1-p}}d\mu}.
\end{equation*}
Then $\psi_{0}\in \mathcal{A}$ and, by (\ref{JJ}), we have
\begin{equation}\label{must1}
J[\psi_{0}]=\left(\int_{\Omega}
u^{\frac{1}{1-p}}d\mu\right)^{-p}\int_{\Omega}
u^{1+\frac{p}{1-p}}d\mu=
\left(\int_{\Omega}
u^{\frac{1}{1-p}}d\mu\right)^{1-p}\geq 1.
\end{equation}
Since $p>1$ then $u$ must satisfy
\begin{equation*}
\int_{\Omega}
u^{\frac{1}{1-p}}d\mu\leq 1.
\end{equation*}
\indent One way to find $\psi_{0}$ is via seeking
equality in the following implication of the reverse H\"{o}lder's inequality:
\begin{equation}\label{fn1}
J[\psi]\geq
\left(\int_{\Omega}\psi\,d\mu\right)^{p}
\left(\int_{\Omega} u^{\frac{1}{1-p}}\,d\mu\right)^{1-p}.
\end{equation}
Equality occurs in (\ref{fn1}) if
$u\psi^{p}=\beta u^{\frac{1}{1-p}},\,$
$\beta$ is a constant. That is
\begin{equation*}
\psi=\beta^{\frac{1}{p}}\, u^{\frac{1}{1-p}}.
\end{equation*}
For $\psi$ to be in $\mathcal{A}$, we must have
\begin{equation*}
\beta=\left(\int_{\Omega} u^{\frac{1}{1-p}}\,d\mu\right)^{-p}
\end{equation*}
and, for this $\beta$, we find $\psi=\psi_{0}$.\\
\indent Another way to identify $\psi_{0}$
is finding a minimizer $\psi_{\text{min}}$ for the nonlinear functional $J$
in $\mathcal{A}$, a simple isoperimetric variational problem:\\
\indent Consider the augmented functional
\begin{equation*}
F(\psi):=\int_{\Omega}
\left(u\psi^{p}-\gamma\psi\right)\,d\mu,\;\;
\psi\in \mathcal{A}.
\end{equation*}
Then, the minimizer $\psi_{\text{min}}$, if it exists, satisfies Euler's equation $p u\psi^{p-1}-\gamma=0$. That is $\psi_{\text{min}}=\left({\gamma}/{p}\right)^{\frac{1}{p-1}}
u^{\frac{1}{1-p}}$. From the restriction
$\int_{\Omega}\psi_{\text{min}}\,d\mu=1$, we can determine
\begin{equation*}
\gamma^{\frac{1}{p-1}}=\frac{p^{\frac{1}{p-1}}}
{\int_{\Omega}\,u^{\frac{1}{1-p}}d\mu}
\end{equation*}
and discover that $\psi_{\text{min}}=\psi_{0}$. Indeed,
if $\int_{\Omega}\,u^{\frac{1}{1-p}}d\mu=1$ then, from the calculation in (\ref{must1}), we get $J[\psi_{\text{min}}]=1$.
\end{proof}
\begin{theorem}\label{main2}
Let $0<p<1$. There does not exist a function $u$
independent of $f$ such that
\begin{equation}\label{pg12}
\left\|
\int_{\Omega} f\,d\mu
\right\|^{p}\,\leq\,\int_{\Omega} u \|
 f \|^{p}\,d\mu
\end{equation}
holds for all $f\in L^{1}(\Omega\rightarrow E)$.
\end{theorem}
\begin{proof}
Assume the contrary. Then, by Lusin's theorem (\cite{fp}, Chapter 2), $\Omega$ possesses
a subset $S$ with $\mu(S)>0$ such that
$u\in L^{\infty}(S\rightarrow \mathbb{R})$.
Pick $\epsilon>0$ so small that $S$ has a subset $S_{\epsilon}$ with $\mu(S_{\epsilon})=\epsilon$.
Let $f= v\,\chi_{S_{\epsilon}}$ so that
$\int_{\Omega} f\,d\mu =\mu\left({S_{\epsilon}}\right)v
=\epsilon\, v$, $\|f\|=\|v\|\,\chi_{S_{\epsilon}}$. From (\ref{pg12}) we have $\int_{S_{\epsilon}} u\,d\mu\geq\epsilon^{p}$. But
$\int_{S_{\epsilon}} u\,d\mu\leq
\mu\left({S_{\epsilon}}\right)\|u\|_{L^{\infty}(S)}=
\epsilon \|u\|_{L^{\infty}(S)}$. Thus
$\|u\|_{L^{\infty}(S)}\geq \epsilon^{p-1}$
which blows up as $\epsilon\rightarrow 0^{+}$ when
$p<1$.
\end{proof}
\section{The discrete inequality}
Let $(X,\|.\|)$ be a normed space
and $(x_{i})_{i\geq1}\subset X$.
Suppose that $\alpha_{i}\in\mathbb{R},$ $i\geq 1$, and $\sum_{i\geq1}\|x_{i}\|<\infty$.
\begin{theorem}
The generalized triangle inequality
\begin{equation}\label{disc1}
\|\sum_{i\geq 1} x_{i}\|^{p}
\leq \sum_{i\geq 1}{{\alpha_{i}} \| x_{i}\|^{p}}
\end{equation}
holds if and only if
\begin{numcases}{}
\alpha_{i}>0,\;\sum_{i\geq 1}
\alpha_{i}^{\frac{1}{1-p}}\leq 1, &$p>1$;\label{suffnec11}\\
\alpha_{i}\geq 1, & $0<p\leq1$.\label{suffnec12}
\end{numcases}
\end{theorem}
\begin{proof}
Let $p>1$. Assume (\ref{suffnec11}) and let $\sum_{i\geq 1} x_{i}\neq 0$. Then, employing the reverse H\"{o}lder's inequality, we obtain
\begin{eqnarray*}
\sum_{i\geq 1}{{\alpha_{i}} \left(\frac{\| x_{i}\|}{\|\sum_{i\geq 1} x_{i}\|}\right)^{p}}
&\geq& \left(\sum_{i\geq 1}{ \frac{\| x_{i}\|}{\|\sum_{i\geq 1} x_{i}\|}}\right)^{p}
\left(\sum_{i\geq 1}{{\alpha^{\frac{1}{1-p}}_{i}}}\right)^{1-p}\\
&\geq& \left(\sum_{i\geq 1}{{\alpha^{\frac{1}{1-p}}_{i}}}\right)^{1-p}
\geq 1
\end{eqnarray*}
by the triangle inequality
\begin{equation}\label{trineqdis}
\|\sum_{i\geq 1} x_{i}\|\leq \sum_{i\geq 1}{\| x_{i}\|}.
\end{equation}
\indent
Suppose (\ref{disc1}) is true. The necessity
of $\alpha_{i}>0$ is clear.
Fix $n$ and let $x_{i}=0$ for all
$i\geq n+1$. Condition (\ref{pgt11}) guarantees
$\sigma_{n}:=\sum_{i=1}^{n}\alpha_{i}^{\frac{1}{1-p}}\leq 1$. Since the sequence $(\sigma_{n})$ is strictly increasing, then
\begin{equation*}
\sum_{i\geq1}\alpha_{i}^{\frac{1}{1-p}}
=\lim_{n\rightarrow \infty} \sigma_{n}\leq 1.
\end{equation*}
The sufficiency of (\ref{suffnec12}) for (\ref{disc1}) is evident from (\ref{trineqdis}) and
the summability of $\left(\|x_{i}\|\right)_{i\geq1}$. Moreover, the validity of (\ref{pleq11}) for
every $n\geq 1$ yields the necessity of (\ref{suffnec12}).
\end{proof}
\begin{theorem}
The reverse generalized triangle inequality
\begin{equation}\label{disc2}
\|\sum_{i\geq 1} x_{i}\|^{p}
\geq \sum_{i\geq 1}{{\alpha_{i}}\| x_{i}\|^{p}}
\end{equation}
holds if and only if $\alpha_{i}<0$ or
\begin{numcases}{\;\exists! j:\alpha_{j}>0,\:
\alpha_{i}<0,\,i\neq j,}
\alpha^{\frac{1}{1-p}}_{j}\leq 1 +\sum_{
i\neq j}|\alpha_{i}|^{\frac{1}{1-p}}, & $p>1$;\label{necsuff21}\\
\alpha_{j}\leq
\inf\{1,\,|\alpha_{i}|,i\neq j\}, & $0<p\leq1$.\label{necsuff22}
\end{numcases}
\end{theorem}
\begin{proof}
Sufficiency of the condition $\alpha_{i}<0$ for
(\ref{disc2}) is obvious
and its necessity results from Theorem \ref{Dadi}.
This leaves us (\ref{necsuff21}) and (\ref{necsuff22}).\\
\indent Assume (\ref{necsuff21}) holds.
Fix $j\geq 1$ and let
$\alpha_{j}>0,$ $\,\alpha_{i}<0,$ $\forall\, i\neq j$. Then, invoking the reverse H\"{o}lder's inequality gives
\begin{align}
\nonumber &
\|\sum_{i\geq 1} x_{i}\|^{p}
+\sum_{i\neq j}{|\alpha_{i}|
\| x_{i}\|^{p}}
\geq\\
&\qquad\qquad\qquad\label{720}
\left(1+\sum_{i\neq j}{|\alpha_{i}|^{\frac{1}{1-p}}}\right)^{1-p}
\left(\|\sum_{i\geq 1} x_{i}\|+
\sum_{i\neq j}{\| x_{i}\|}\right)^{p}.
\end{align}
Notice that
\begin{equation}\label{721}
\|\sum_{i\geq 1} x_{i}\|+
\sum_{i\neq j}{\| x_{i}\|}\geq
\| x_{j}\|.
\end{equation}
By (\ref{necsuff21}) and (\ref{721}), the inequality (\ref{720}) simplifies to
\begin{equation*}
 \|\sum_{i\geq 1} x_{i}\|^{p}
+\sum_{i\neq j}{|\alpha_{i}|
\| x_{i}\|^{p}}\geq \alpha_{j}\| x_{j}\|^{p}
\end{equation*}
from which follows (\ref{disc2}). \\
\indent Conversely, assume (\ref{disc2}).
Let $ \rho_{n}=\sum_{\substack{i=1,...,n,\\
i\neq j}}|\alpha_{i}|^{\frac{1}{1-p}}$.
Condition (\ref{pgt1}) asserts
that $\alpha_{j}\leq 1+\rho_{n}$ for all $n\geq 2$.
The limit of $\rho_{n}$,
even though can be infinite, exists by
$(\rho_{n})$ being strictly increasing.
Thus
\begin{equation*}
\alpha_{j}\leq 1+\lim_{n\rightarrow \infty}\rho_{n}
=1+\sum_{i\neq j}|\alpha_{i}|^{\frac{1}{1-p}}.
\end{equation*}
Observe that (\ref{720}) tolerates
the case $\sum_{i\neq j}|\alpha_{i}|^{\frac{1}{1-p}}=\infty$
because $1-p<0$.\\
\indent Let $0<p\leq 1$. The necessity of (\ref{necsuff22}) for (\ref{disc2}) follows easily from (\ref{pleq1}). Next, assume (\ref{necsuff22}). Then
\begin{equation}\label{cof1}
\frac{1}{\alpha_{j}}\geq1,\;
\frac{|\alpha_{i}|}{\alpha_{j}}\geq 1,\,\forall\, i\neq j.
\end{equation}
The case $x_{j}=0$ returns us to (\ref{disc2}) with
$\alpha_{i}<0,\;\forall\, i\geq 1$. So
let $x_{j}\neq 0$. Now since
\begin{equation*}
\min_{\substack{\\ \xi_{i}\geq 0,\,\forall\, i\geq 1,\\ \sum_{i\geq 1}\xi_{i}\geq 1}}\sum_{i\geq 1}\xi^{p}_{i}= 1.
\end{equation*}
Then, exploiting (\ref{721}), we obtain
\begin{equation}\label{disc2p11}
\left(\frac{\|\sum_{i\geq 1} x_{i}\|}{\|x_{j}\|}\right)^{p}
+\sum_{i\neq j}\left(\frac{\| x_{i}\|}{\|x_{j}\|}\right)^{p}
\geq 1.
\end{equation}
Finally, introducing (\ref{cof1}) to (\ref{disc2p11}), we deduce
\begin{equation*}
\frac{1}{\alpha_{j}}\left(\frac{\|\sum_{i\geq 1} x_{i}\|}{\|x_{j}\|}\right)^{p}
+\sum_{i\neq j}\frac{|\alpha_{i}|}{\alpha_{j}}
\left(\frac{\| x_{i}\|}{\|x_{j}\|}\right)^{p}
\geq 1
\end{equation*}
which implies (\ref{disc2}).
\end{proof}
\section*{Acknowledgement}
The author thanks the referees for
their careful reading of the manuscript.

\end{document}